\documentclass[12pt]{article}
\usepackage{amsfonts}
\newcommand{\C}{{\bf C}}

\newcommand{\BP}{{\bf P}}
\newcommand{\Q}{{\bf Q}}
\newcommand{\nind}{\noindent}
\newcommand{\rar}{\rightarrow}
\newcommand{\olin}{\overline}
\newcommand{\sms}{\setminus}
\newcommand{\A}{{\bf A}}
\newcommand{\tild}{\widetilde}
\newcommand{\lkd}{\overline\kappa}
\begin{document}
\title{Non-singular affine surfaces with self-maps}
\author{R.V. Gurjar and D.-Q. Zhang}
\maketitle
\begin{abstract}
We classify surjective self-maps (of degree at least two) of affine
surfaces according to the log Kodaira dimension.
\end{abstract}

\par \vskip 1pc \noindent
{\bf \S 1. Introduction.}\\

In this paper we are interested in the following question.\\

\nind
{\bf Question.} Classify all smooth affine
surfaces $X$/${\C}$ which admit a proper morphism $f:X\rar X$ with degree
$f>1$.\\

In \cite{Fu} and \cite{Na}, a classification of smooth projective surfaces
with a self-map of degree $>1$ has been given. This paper is inspired by
their results. The case when $X$ is singular appears to be quite hard so
we restrict ourselves to the smooth case. Similarly, if $f$ is not a
proper morphism then again the problem is difficult. For example, we do
not even know if there is an {\'e}tale map of degree $>1$ from $\C^2$ to
itself. This is the famous Jacobian Problem. If $S$ is any $\Q$-homology
plane with $\lkd(S)=-\infty$ then $S$ admits an algebraic action of
the
additive group $G_a$ (\cite{MM}). Hence the automorphism group of such a
surface is
infinite. However, the problem of constructing a proper self-map of degree
$>1$ for $S$ is quite non-trivial. Our main result can be stated
as follows.\\

\nind
{\bf Theorem.} {\it There is a complete classification of smooth complex
affine surfaces $X$ which admit a proper self-morphism of degree $>1$,
if either the logarithmic Kodaira dimension
$\lkd(X)\geq 0$ or the topological fundamental group $\pi_1(X)$ is infinite.
\par
More precisely, any such $X$ is isomorphic to a quotient of the form $(\Delta\times \A^1)/G$ or
$(\Delta\times\C^*)/G$ where $\Delta$ is a smooth curve and $G$ is a finite group acting freely on
$\Delta\times\A^1$ or $\Delta\times\C^*$ respectively.}\\

\par
As a consequence of the proof, we have:

\par \vskip 1pc \noindent
{\bf Corollary.} Suppose that $X$ is an affine surface with a proper morphism $X \rightarrow X$
of degree $> 1$. Then we have:
\begin{itemize}
\item[(1)]
If $\lkd(X) \ge 0$, then $\lkd(X) = 0, 1$
and $X \cong (\Delta\times\C^*)/G$ where $\Delta$ is a smooth affine curve and
$G$ is a finite group acting freely on $\Delta\times\C^*$.

\item[(2)]
Suppose that $\lkd(X) = -\infty$ and let
$\varphi : X \rightarrow B$ be an ${\bold A}^1$-fibration
(for the existence see [15], Chapter I, \S 3).

\item[(2a)]
If $\lkd(B) = -\infty$, then the topological fundamental group $\pi_1(X)$
is finite.

\item[(2b)]
If $\lkd(B) = 0$, then $B  \cong {\bold C}^*$ and $X \cong {\bold A}^1 \times {\bold C}^*$.

\item[(2c)]
If $\lkd(B) = 1$ then every fibre of $\varphi$ is reduced and irreducible
and $X \cong (\Delta \times {\bold A}^1)/({\bold Z}/(m))$ ($m \ge 1$)
where $\Delta$ is a smooth affine curve
and ${\bold Z}/(m)$ acts freely on $\Delta \times {\bold A}^1$.
\end{itemize}
\par \vskip 1pc

\nind
{\it Remark.}
\par
(1) There are easy examples of $X$
where $\lkd(X) \ge 0$ or $\pi_1(X)$ is infinite, but
$X$ has no proper self map of degree $> 1$. For example, let $C,D$ be smooth
irreducible projective curves of genus $m,n\geq 2$ which admit no non-trivial automorphisms.
Assume that $m\neq n$ and let $C',D'$ be non-empty affine open subsets of $C,D$ respectively. Then
$C',D'$ have no proper self-maps of degree $>1$ and hence $X:=C'\times D'$ has no proper self-map of
degree $>1$. Clearly, $\lkd(X)=2$ and $\pi_1(X)$ is infinite.\\
It is an interesting problem to determine
those quotients $(\Delta\times \A^1)/G$ and
$(\Delta\times\C^*)/G$ which have no proper self maps of degree $> 1$.

\par
(2) So far the authors have not been able to find any example
(apart from $\C^2$) of a smooth
affine surface $X$ with $\lkd(X)=-\infty$ and $\pi_1(X)$ finite
such that $X$ has a proper self-map of degree $>1$.

\par
(3) Our proof of the theorem uses the classification theory of open algebraic
surfaces
developed by S. Iitaka, Y. Kawamata, T. Fujita, M. Miyanishi,
T. Sugie and other Japanese mathematicians. We also use topological
arguments in an essential way.\\

\par \vskip 1pc \noindent
{\bf Acknowledgement.}
The authors are thankful to Peter Russell for suggesting the argument in Lemma 3.4A from
C.H. Fieseler's paper \cite{Fi}.
The authors also like to thank the referees for constructive suggestions
which improve the paper.
\\

\nind
{\bf \S 2. Preliminaries.}\\

We will only deal with complex algebraic varieties. By a curve
(resp. surface) we mean
an irreducible, quasi-projective curve (resp. an irreducible,
quasi-projective surface). By a {\it component} of a variety $Z$ we mean
an
irreducible component of $Z$. Let $Z$ be a smooth surface.
By a ${\BP}^1$-fibration on $Z$ we mean a morphism $Z\rar B$ onto a smooth
algebraic curve whose general fiber is isomorphic to ${\BP}^1$. Similarly
an $\A^1$-fibration and a $\C^*$-fibration on $Z$ can be defined. Here,
$\C^*$ denotes $\BP^1-~\{two~points\}$.\\
Given a fibration $\varphi:X\rar B$ from a smooth surface $X$ onto a smooth curve $B$ let $F:=\sum m_iF_i$
be a scheme-theoretic fiber of $\varphi$ where $F_i$ are the irreducible components of $F$. The greatest
common divisor of the integers $m_i$ is called the {\it multiplicity of $F$}. If the multiplicity is
$>1$ then we say that $F$ is a {\it multiple fiber}.\\
A smooth projective irreducible rational curve $C$ on a smooth surface $Z$ with $C^2=n$ is called an
$(n)$-curve. The topological Euler-Poincar{\'e} characteristic of a variety
$Z$ is denoted by $\chi(Z)$.\\
Given a smooth surface $Z$ there is an open embedding $Z\subset W$
such that $W$ is a smooth projective surface and $D:=W\sms Z$ is a divisor with
simple normal crossings. If any $(-1)$-curve in $D$ meets at least three
other components of $D$ then we say that $W$ is a {\it good}
compactification of $Z$.

For a smooth, irreducible, quasi-projective variety $X$,
we denote by $\lkd(X)$ the logarithmic Kodaira dimension.
So when $X$ is compact, then $\lkd(X)$ is just the usual one,
and when $X$ is non-compact, we let $\overline{X}$ be a compactification
with $D = \overline{X} \setminus X$ a simple normal crossing divisor
and then $\lkd(X) = \lkd(\overline{X}, K_{\overline X} + D)$
as defined by Iitaka, which is independent of the choice of the compactification.
\\

\nind
{\bf \S 3. Proof of the Theorem}.\\

Let $X$ denote a smooth affine surface which admits a proper morphism
$f:X\rar X$ with degree $f>1$. Our proof of the classification splits into
many cases (see \S 3.1 $\sim$ \S 3.7 below).\\

We will begin
with the easier case when $\lkd(X)\geq 0$. Then by a basic result
of Iitaka \cite{I} the map $f$ is {\'e}tale. Since $f$ is also proper
by assumption it follows by covering space theory that
$\chi(X)=deg~f\cdot\chi(X)$. Hence $\chi(X)=0$.\\

\nind
{\bf 3.1.} Suppose that $\lkd(X)=2$.

\par
We claim that this case cannot occur. As a corollary of an
inequality of Miyaoka-Yau type by R. Kobayashi, S. Nakamura and F. Sakai
it follows that in this case $\chi(X)>0$ (\cite{MT}). Hence $f$ cannot
exist in this case.\\

\nind
{\bf 3.2.} Suppose that $\lkd(X)=1$.
\par
Since $X$ is affine, by a basic result due to Kawamata there is a
$\C^*$-fibration $\varphi:X\rar B$ (\cite{M}, Chapter II, \S 2). Using
Suzuki's formula we
can calculate $\chi(X)$ in terms of $\chi(B)$ and $\chi(F_i)$, where
$F_i$ are the singular fibers of $\varphi$ (\cite{Su} and \cite{Z}). Since
$\chi(X)=0$, it follows from this formula that every fiber of $\varphi$ is
isomorphic to $\C^*$, if taken with reduced structure. Let
$m_1F_1,m_2F_2,\ldots,m_rF_r$ be all the multiple fibers with multiplicities $m_i$ of $\varphi$.\\

\par
Suppose first that $r=0$. The $\C^*$-fibration $\varphi$ may not be
Zariski locally-trivial but there is a $2$-sheeted {\'{e}}tale covering $\Delta\rar B$ such that
the
fiber product $\Delta\times_B X$ is a Zariski locally-trivial fibration.
Then it is easy to see that $\lkd(X)=1$
implies that
$\lkd(B)=1$ \cite{K1}. Hence either $B$ is a non-rational curve or a rational
curve with at least three
places at infinity. It follows by L{\"u}roth's theorem that $f$ maps any
fiber of $\varphi$ onto another fiber of $\varphi$. This induces an
{\'e}tale, proper self-map $f_0:B\rar B$. Now $\chi(B)<0$, since
$\lkd(B)=1$. Hence $f_0$ is an automorphism which must have finite
order since $B$ is of general type. Hence $f^N$ induces an identity on
$B$ for some $N>0$. Let $g:=f^N$. In view of Claim 3.2a below we will assume that $\varphi$ is itself
locally trivial. We can find an open embedding $X\subset W$ such that $W$
has a $\BP^1$-fibration $\olin\varphi:W\rar B$ and $W\sms X$ is a disjoint
union of two cross-sections $S_1,S_2$ of $\olin\varphi$. The map $g$
extends to $W$ mapping every fiber of $\olin\varphi$ to itself. By
considering $g^2$, if necessary, we can assume that $S_1,S_2$ are
pointwise fixed.\\

\noindent
{\bf Claim 3.2a.} There is a finite {\'e}tale cover $\Delta'\rar B$ such that
$X\times_B\Delta'\cong\Delta'\times\C^*$.\\

\par
{\it Proof of the claim.}
We may assume that $\Delta = B$.
Restricted to each fiber of $\varphi$, the map $g$ has the form $t\rar
at^m$, where $a$ is a non-zero complex number and $m$ is an integer $\geq
2$ (since deg $g>1$). We can cover $B$ by open subsets $U_1,U_2$ such that
$\varphi$ is trivial on each $U_i$. Let $\varphi$ be obtained by patching
$U_1\times\C^*,~U_2\times\C^*$ in $(U_1\cap U_2)\times\C^*$ by $(z,t)\sim
(z,\eta(z)t)$, where $\eta$ is a regular invertible function on $U_1\cap
U_2$. The map $g$ can be described on $U_i\times\C^*$ by $g(z,t)=(z,\alpha_i
(z)t^m)$, where $\alpha_i$ is a unit on $U_i$. Patching gives the relation
$\alpha_1\eta = \alpha_2 \eta^m$, i.e. $\alpha_1/\alpha_2=\eta^{m-1}$. Hence
the $\C^*$-bundle $\varphi$ with transition function $\eta$ is torsion of
order $m-1$. It follows that there is an {\'e}tale cover $\Delta'\rar B$
of degree $m-1$ with the required property. This proves the claim.\\

\par
Now we know that $X\cong (\Delta'\times\C^*)/H$, the action of
$H$ being fixed point free. Here $H$ is a finite group of
degree $m$ or $2m$.
\\

\par
Next assume that $r>0$. By the solution of Fenchel's conjecture due to Fox
and Bundgaard-Nielson (see \cite{BN}, \cite{Ch}, \cite{Fo}),
there exists a Galois covering $\tilde
B\rar B$, ramified precisely over $\varphi(F_i)$ with ramification index
$m_i$ for $i=1,2,\ldots,r$. Then the normalization $Y = \olin{X\times_B \tilde B}$
of the fiber
product $X\times_B\tilde B$ is an {\'e}tale cover of $X$. The
induced $\C^*$-fibration $Y\rar \tilde B$ has no singular fibers.
Hence
$Y\cong (C \times \C^*)/H$ for some {\'e}tale cover
$C\rar \tilde B$ of degree equal to $|H|$. By going to a further covering $C'$ of $C$,
if necessary, we can assume that $C'\rar B$ is finite Galois and
$\olin{X\times_B C'}$ is an {\'e}tale Galois covering of $X$. Thus $X$ is
a quotient of a surface of the form $C'\times\C^*$ by a fixed point free
action of a finite group.\\

\par
This completes the description of $X$ when $\lkd(X)=1$.\\

\nind
{\bf 3.3.} Assume now that $\lkd(X)=0$.\\

Since $\chi(X)=0$, we see that $b_1(X)\neq 0$. Let ${\cal A}(X)$ be the
quasi-Albanese variety of $X$ and let $\alpha:X\rar{\cal A}(X)$ be the
quasi-Albanese map (\cite{K2}). There is a short exact sequence:
$$(0)\rar (\C^*)^l\rar{\cal A}(X)\rar Alb(V)\rar (0)$$
where $V$ is a smooth projective compactification of $X$ and $Alb(V)$ the
Albanese variety of $V$. Since $\lkd(X)=0$, by a result of Kawamata
$\alpha$ is a dominant map with a connected general fiber (\cite{K2}).\\

\par
We assert that if the image of $X\rar{\cal A}(X)$ is a curve then
$\lkd F = 0$ with $F$ a general fibre of $\alpha$, whence
$F \cong \C^*$ because $X$ is affine.
Indeed, the easy addition for log Kodaira dimension implies that
${\lkd}(F) \ge 0$. On the other hand, Kawamata also proved
${\lkd}(X) \ge {\lkd}(F) + {\lkd}({\cal A}(X))$.
So to show the assertion, we have only to show that
${\lkd}({\cal A}(X)) \ge 0$. To see this, note that ${\cal A}(X)$
is either isomorphic to the Abelian variety $Alb(V)$, or to ${\bold C}^*$
whence ${\lkd}({\cal A}(X)) \ge 0$. This also proves the assertion.

If the image of
$\alpha$ is a surface then we consider the composite map $X\rar{\cal
A}(X)\rar Alb(V)$. If the image of this composite map is a curve then
again $X$ has a $\C^*$-fibration.

\par
In this case and in the case where the image of $X\rar{\cal A}(X)$ is a curve,
$X$ has a ${\bold C}^*$-fibration. Then
we argue as in the
previous case. From $\chi(X)=0$ we see that the $\C^*$-fibration has
at most multiple $\C^*$'s as singular fibers. We can again conclude that $X$ is
a quotient of $\Delta\times\C^*$ by a finite fixed point-free
automorphism group for a suitable curve $\Delta$.\\

\par
Assume next that the image of $X\rar Alb(V)$ is a surface.
We shall show that this case does not occur.
Indeed, note that the image of
$\alpha$ is also a surface and by \cite{K2} $\alpha$ is birational.
Thus $2 = \dim X \ge \dim {\cal A}(X) \ge \dim Alb(V) \ge 2$, whence in the above exact sequence,
one has $l = 0$. So ${\cal A}(X) = Alb(X)$.
Thus the map $X\rar Alb(V)$ is birational. There is a curve
$C$ in $X$ (which may be reducible or empty) such that $C$ maps to
finitely many (smooth) points in $Alb(V)$ such that $X\sms C$ is an open
subset of $Alb(V)$. Then $\lkd(X\sms C)=\lkd(X)=0$. But
$X\sms C$ is an affine subset of the abelian variety $Alb(V)$ so its complement
in $Alb(V)$ is an ample divisor (see \cite{Ha}).
Hence we see easily that $\lkd(X\sms C)=2$, a contradiction.\\

\par
Finally, if the map $X\rar Alb(X)$ is trivial then there is a non-constant
map $X\rar (\C^*)^l$, and hence a map $X\rar\C^*$ which must then be a
$\C^*$-fibration by Kawamata's inequality in \cite{K1} and by ${\lkd}(X) = 0$.
Thus, when $\lkd(X)=0$, as in the previous case,
$X\cong Y/G$ for some smooth affine surface $Y$ which is isomorphic to
$\Delta\times\C^*$ for some curve $\Delta$ and $G$ acting fixed point
freely on $Y$.\\

\nind
{\bf 3.4.} Now assume that $\lkd(X)=-\infty$.

\par
Note that $X$ has a compactification so that the complement of $X$
supports an ample divisor (and hence connected) by \cite{Ha}.
By a basic result due to Fujita, Miyanishi and Sugie (\cite{M}, Chapter
I, \S 3), there is
an $\A^1$-fibration $\varphi:X\rar B$.\\

\par \vskip 1pc \noindent
{\bf Proposition 3.4.} Suppose that $\lkd(X) = -\infty$
and $\varphi : X \rightarrow B$ is an ${\bold A}^1$-fibration
such that $\lkd(B) \ge 0$. Then every fibre of $\varphi$ is reduced and
irreducible.

\par \vskip 1pc
We now prove the proposition. So suppose that $\lkd(B) \ge 0$.
Note that $f$ permutes the fibers of
$\varphi$.

\par
We first treat the case when some fiber $F_0$ of $\varphi$ is not
irreducible, though a general fibre is a reduced curve isomorphic to $\A^1$.\\

\par
Then $f^{-1}(F_0)$ is also not irreducible since $\varphi$ is
surjective. As before there is an induced map $f_0:B\rar B$. By
considering $f^N$ for suitable $N$ we can assume that $f_0$ is an identity
map and every irreducible component of $F_0$ is mapped to itself by $f$.

In this case we prove a general result which will be useful
in later arguments as well.\\

\nind
{\bf Lemma 3.4A.} {\it There is no proper self-map $f:X\rar X$ which maps any
irreducible component of any fiber of $\varphi$ to itself.}\\

\nind
{\it Proof.} For an irreducible component $C$ of $F_0$, let $f^*(C)=aC$.
Since the
induced map on $B$ is identity we easily deduce that $a=1$. Since $X$
is affine the fiber $F_0$ is a disjoint union of curves isomorphic to
$\A^1$. It follows that no irreducible component of any fiber of $\varphi$
is ramified for the map $f$. Let $X_0$ be obtained from $X$ as follows. For any
reducible fiber of
$\varphi$ other than $F_0$ we remove all but one
irreducible
components. If $F_0$ is non-reduced then we remove all the irreducible components of $F_0$
except for one non-reduced irreducible component. If $F_0$ is reduced then we
remove all but two irreducible components of $F_0$. Using the solution of Fenchel's conjecture
by Bundagaard-Nielsen, Fox
\cite{BN}, \cite{Fo} and noting
that $B \ne {\bold P}^1$,
we can construct a finite Galois ramified cover $\tild B\rar B$
with prescribed ramification such that the normalization of the fiber
product $X_0\times_B\tild B$, say $\tild X$, is a finite {\'e}tale cover of $X_0$.
By the universal property of fiber products it follows that $f$ lifts to a finite self-map
$\tild f: \tild X\rar \tild X$.
The induced $\A^1$-fibration on $\tild X$ has all fibers reduced and at
least one fiber not irreducible. From this observation it now suffices to
assume that $X_0=\tild X$.\\

Let $C_1,C_2$ be distinct irreducible components of $F_0$. Denote by $X_0$
the affine surface obtained from $X$ by removing all but one
irreducible
components of all the fibers of $\varphi$, other than $F_0$, and all the
irreducible components of $F_0$, other than $C_1,C_2$. Then there is an
induced proper self-map $f_0:X_0\rar X_0$ which maps every irreducible
fiber-component to itself. Every fiber of the $\A^1$-fibration is
reduced. Let $S$ be the affine scheme over the power series ring in one variable
$k[[t]]$ obtained from $X_0$ with an $\A^1$-fibration over $k[[t]]$ whose special fiber is
reduced with irreducible components $C_1,C_2$. To
complete the proof of Lemma 3.4A, it suffices to prove that there is no finite
self-map $f:S\rar S$ of degree $>1$.\\

For this we use an idea from Fieseler's paper (\cite{Fi}, \S 1). Let
$U:=Spec~k[[t]]$, $~U^* := U-\{(t)\}$. Since the special fiber is reduced
it is clear that $S-C_1\cong U\times\A^1,~S-C_2\cong U\times\A^1$. Now $S$ is obtained by
patching $S-C_1,~S-C_2$ along $U^*\times\A^1$. Let the patching map $p:(S-C_1)-C_2\rar
(S-C_2)-C_1$ be given by $p(t,y)=(t,\alpha(t)y+\beta(t))$, where $\alpha,\beta\in k((t))$ and
$\alpha\neq 0$. The self-map $f$ restricted to $S-C_1$ is given by
$f(t,y)=(t,a_0(t)y^m+a_1(t)y^{m-1}+\cdots+a_m(t))$. Each $a_i\in k[[t]]$.
Since $f$ has degree
$m$
restricted to each fiber
and also restricted to $C_1,C_2$, our $a_0$ is a unit. Similarly, $f$ restricted to $S-C_2$ is given
by $f(t,y)=(t,b_0y^m+\cdots+b_m)$. Using the patching we have $p\circ f=f\circ p$. This means
$$\beta + \alpha (a_0 y^m + a_1 y^{m-1} + \dots + a_m) =
b_0 (\alpha y + \beta)^m + b_1 (\alpha y + \beta)^{m-1} + \dots + b_m.$$
Equating the first two leading coefficients, we get:
$$\alpha a_0 = b_0 \alpha^m,$$
$$\alpha a_1 = b_0 m \alpha^{m-1} \beta + b_1 \alpha^{m-1}.$$
Solving them, we obtain:
$$\alpha^{m-1} = a_0/b_0,$$
$$(ma_0) \beta = \alpha a_1 - a_0 b_1/b_0.$$
Since $a_0$ and $b_0$ are units,
it follows that $\alpha,\beta\in k[[t]]$. Since the special fiber is a
disjoint union of $C_1,C_2$ the map $p$ cannot extend to a morphism $S-C_1\rar S-C_2$.
Therefore $\alpha,\beta$ cannot be both in $k[[t]]$. This contradiction completes
the proof of Lemma 3.4A.\\

\par
Now assume that every fiber of $\varphi$ is irreducible.\\

If further some fiber is non-reduced then by going to a suitable ramified
cover $\Delta\rar B$ we see that the normalization of the fiber product
$Y:=\olin{X\times_B\Delta}$ is an {\'e}tale cover of $X$ and has an
$\A^1$-fibration over $\Delta$ with all
reduced fibers and at least one fiber which is not irreducible. By
construction $f$ lifts to a proper self-map $f':Y\rar Y$ of degree $>1$.
By the previous case such a map cannot exist.
This completes the proof of the proposition.\\

\par
Now assume that every fiber of $\varphi$ is reduced and irreducible.
We consider the cases $\lkd(B) = 1, 0$ separately. \\

\par
Consider first the case $\lkd(B)=1$.\\

The induced map $f_0:B\rar B$ is an {\'e}tale, finite
map, hence an automorphism of finite order.
By taking $f^N$ for suitable $N\geq 1$ we will assume that $f_0$ is
identity.

In this case we will prove that there is a finite {\'e}tale cover $\tild B\rar B$ such that the fiber
product $\tild B\times_B X$ is a trivial $\A^1$-bundle over $\tild B$.

\nind
As above we can embed $X$ into a smooth surface $V$ with a
$\BP^1$-fibration $\Phi:V\rar B$ extending $\varphi$ such that $V-X$ is a cross-section
$S$ of $\Phi$.
It is easy to see that $f$ extends
to a self-map $V\rar V$ which we write again by $f$ for simplicity. Now $f$ maps any point in $S$ to
itself. We argue as in the proof of
claim 3.2a. Let $B=U_1\cup U_2$ be an open cover such that the $\A^1$-fibration on each $U_i$ is
trivial. Let the patching be given on $U_1\cap U_2$ be $(z,t)\sim (z,\eta(z)t+\xi(z))$, where
$\eta$ is a
nowhere vanishing regular function on $U_1\cap U_2$ and $\xi$ is regular on $U_1\cap U_2$. Let $f$
over $U_1$ be given by $f(z,t)=(z,a_0(z)t^d+a_1(z)t^{d-1}+\cdots+a_d(z))$ and over $U_2$ by
$f(z,t)=(z,b_0(z)t^d+b_1(z)t^{d-1}+\cdots + b_d(z))$. Here $a_i,b_j$ are regular functions on $U_1, U_2$
respectively
and $a_0,b_0$ are nowhere zero on $U_1,U_2$ respectively. Using patching we get on $U_1\cap U_2$ the
following $$\eta(z)(a_0t^d+a_1t^{d-1}+\cdots+a_d)+\xi=b_0(\eta t+\xi)^d+b_1(\eta
t+\xi)^{d-1}+\cdots+b_d.$$
Comparing the coefficients of $t^d$ on both sides we get $\eta a_0=b_0\eta^d$. Hence
$\eta^{d-1}=a_0/b_0$, showing that there is a torsion line bundle on $B$ whose order divides $d-1$.
Consider the {\'e}tale cover $\tild{U_1}$ of $U_1$ obtained by adjoining $\tild{a_0}:=a_0^{1/(d-1)}$ to
the
coordinate
ring of $U_1$ and similarly let $\tild{U_2}$ be obtained by adjoining $\tild{b_0}:=b_0^{1/(d-1)}$ to the
coordinate ring of $U_2$. Since $\tild{a_0}/\tild{b_0}=\eta$ these patch to give an {\'e}tale cover
$\tild
B$
of $B$. The self-map $f$ extends to a proper self-map $\tild f:\tild X\rar \tild X$, where $\tild
X=X\times_B\tild
B$.
We will show that the pull-back $\A^1$-fibration $\tild X$ is a trivial $\A^1$-bundle over
$\tild B$. It is easy to see that the patching on $\tild{U_1}\cap\tild{U_2}$ can be assumed to be of the
form
$\tild p(z,t)=(z,t+\tild\xi)$. Writing the self-map $\tild f$ on $\tild{U_i}$ as above we can assume by
our
construction that the coefficients $\tild {a_0}=\tild {b_0}$. Comparing the
coefficients of $t^{d-1}$ we get $\tild{a_1}=d\tild{b_0}\tild\xi+\tild{b_1}$. Hence
$\tild\xi=(\tild{a_1}-\tild{b_1})/d\tild{b_0}$. Since $\tild{a_0}=\tild{b_0}$, the function
$\tild{b_0}$ is a nowhere vanishing function on whole of $\tild B$. By changing the notation,
$\tild\xi$ is a difference
of two regular functions $\tild{a_1},~\tild{b_1}$ on $U_1,U_2$ respectively. From this it is easy to
deduce that the $\A^1$- fibration on $\tild
B$ is trivial.

\par
Hence we have shown that $X$ is a quotient of a product $\tild B\times\A^1$ by a finite cyclic group
acting fixed point freely.\\

\par
Consider next the case $\lkd(B) = 0$.
If $B\cong\C^*$ then $X\cong {\bold A}^1\times\C^*$. In this case $X$ clearly has
proper self-map of arbitrary degree.\\

Suppose that $B$ is an elliptic curve.
Clearly $f$ permutes the fibers of $\varphi$, thus inducing a self-map $f_0:B\rar B$.
We can choose a relatively minimal ruled surface $V$ as a compactification of
$X$ so that $V \setminus X$ equals the cross-section $S$ at "infinity" of
the unique ruling $V \rightarrow B$. One sees easily that
$f : X \rightarrow X$ extends to a self map of $V$, also denoted by
$f : V \rightarrow V$ (here the irrationality of the base $B$
is essentially used). We shall show that there is no such $f$
of degree $> 1$.

\par \vskip 1pc \noindent
{\bf Lemma 3.4B.} Suppose that ${\lkd}(X) = -\infty$
and the base curve $B$ of the ${\bold A}^1$-fibration
$\varphi : X \rightarrow B$ (each fibre of which is an irreducible and reduced ${\bold A}^1$)
is a nonsingular elliptic curve. Then there is no self map $f : X \rightarrow X$
of degree $> 1$.

\par \vskip 1pc \noindent
{\it Proof.}
Suppose the contrary that $f$ is a self map of $X$ of degree $> 1$.
We use the notation above: $f : V \rightarrow V$, $S = V \setminus X$, etc.
Write $f^*S = eS$ with $e$ the ramification index.
Let $e'$ be the index of function fields extension $|{\bold C}(S) : {\bold C}(f(S))|$.
Then $deg(f) = e e'$. Since
$e e' S^2 = (f^*S)^2 = (eS)^2 = e^2 S^2$ and since $S$
is ample (so $S^2 > 0$),
we have $e = e'$.

We assert that $f : X \rightarrow X$ is etale. Suppose the contrary that
this $f$ is ramified along a divisor $D^0$ and let $D$ be the closure in $V$
(the purity of branch locus over regular ring, is used).
The ampleness of $S$ implies that $S \cap D$ contains a point $P$.
Thus $f | S : S \rightarrow S$ is ramified at $P$.
This is impossible because $S$ is an elliptic curve.
So the assertion is true and $f : X \rightarrow X$ is etale.
In particular, for every fibre $F$ on $X$, the map $f: F \rightarrow f(F)$
is etale and hence an isomorphism because $F \cong {\bold A}^1$ is simply connected.
Since $S$ is ramified the point at infinity on any $F$ is ramified. This
is a contradiction. Hence deg $f=1$. This proves the lemma.\\

\par \vskip 1pc
Now we are left with the case ${\lkd(B)} = -\infty$. Then $B=\A^1$ or $\BP^1$.\\

\nind
{\bf Lemma 3.5.} It is impossible that $B=\A^1$
and the fundamental group $\pi_1(X)$ is infinite.\\

{\it Proof.} Suppose the contrary that $B=\A^1$
and $\pi_1(X)$ is infinite.
We claim that the first Betti number $b_1(X)=0$. This follows from the exact sequence for homology
groups with rational coefficients (\cite{Su}):
$$H_1(\A^1)\rar H_1(X)\rar H_1(B)\rar (0).$$
It follows that $\chi(X)>0$. Since $\pi_1(X)$ is infinite we see that $\varphi$ has at least two
multiple fibers (otherwise, $\pi_1(X)$ is finite cyclic). Let $\Delta\rar B$ be a suitable ramified
Galois covering such that $Y:=\overline{X\times_B\Delta}\rar X$ is {\'e}tale and $\psi:Y\rar\Delta$ is
an $\A^1$-fibration with all reduced fibers. Now $\lkd(\Delta)\geq 0$.\\

\nind
{\it Claim.} $f$ extends to a finite map $g:Y\rar Y$.\\

For this we will show that the induced homomorphism
$f_*:\pi_1(X)\rar\pi_1(X)$ is an
isomorphism. To see this, first we know that this image
always has finite
index. This observation goes back to Serre (cf. \cite{N},
Lemma 1.5). Let
$Z\rar X$ be the finite {\'e}tale cover such that
$\pi_1(Z)$ is equal to
the above subgroup of finite index. By covering space
theory, $f$ extends
to a morphism $f':X\rar Z$. Let $d$ be the index of
$\pi_1(Z)$ in
$\pi_1(X)$. Then $\chi(Z)=d\chi(X)$. Since $X$ dominates
$Z$, by Lemma 1.5
in \cite{N} we know that $0=b_1(X)=b_1(Z)$. Since $f'$ is
also a finite
map it is known that $b_2(X)\geq b_2(Z)$ \cite{Gi}. Now
$\chi(Z)=d\chi(X)\geq d \chi(Z)$, we have a contradiction if
$d\geq 2$. This
proves the assertion that $f_*$ is onto.\\
We claim that, in fact, $f_*$ is an isomorphism.\\
To see this, let $m_1F_1,m_2F_2,\ldots,m_rF_r$ be the
multiple fibers of
$\varphi$. A slight extension of Lemma 1.5 in \cite{N}
shows that we have
an isomorphism (using the fact that a general fiber of
$\varphi$
is $\A^1$)
$$\pi_1(X)\cong
<e_1,e_2,\ldots,e_r|e_1^{m_1}=1=\cdots=e_r^{m_r}>$$
This group is known to be residually finite, i.e. the intersection of all its subgroups of finite index
is trivial (\cite{Gr}). By a result of Malcev any finitely generated residually finite group $G$ has
the property
that any surjective
homomorphism $G\rar G$ is an isomorphism (\cite{Ma}). Hence the
surjection $f_*$ is an isomorphism.

\par
Now by covering space theory $f$ extends to a finite self map $g:Y\rar Y$. This proves the claim. \\

\par
At least two fibers of $\psi$ (over $b_1, b_2 \in \Delta$ say) are
not irreducible.
Note that $g$ induces an automorphism $g | \Delta$ on the base $\Delta$.
A power $h = g^s$ induces $h | \Delta$ fixing $b_1, b_2$.
If the affine curve $\Delta$ has $\lkd(\Delta) \ge 0$, we see that
a higher power $g^{sn}$ induces the identity :$g^{sn} | \Delta = id$;
if $\lkd(\Delta) = -\infty$, then $\Delta = {\bold A}^1 \subset {\bold P}^1 =
\Delta \cup \{\infty\}$ and $g^s$ induces an automorphism of ${\bold P}^1$
fixing $b_1, b_2, \infty$ (so $g^s | {\bold P}^1 = id$), whence $g^s$ stabilizes every fibre of $\varphi.$
Now we can apply
the argument in the proof of Lemma 3.4A and conclude that $g^{sn}$ cannot exist.
Hence $f$ also does not exist. This proves the lemma.\\

\nind
{\bf Lemma 3.6.} It is impossible that $B=\BP^1$
and the fundamental group $\pi_1(X)$ is infinite.\\

The proof in this case is similar to the above case. Indeed, if $\pi_1(X)$ is
infinite, then $\varphi$ has at least three multiple fibers and in the case
of three multiple fibers the
multiplicities do not form a Platonic triple. Now the rest of the argument
is very similar to the previous case.\\
This completes the proof of the Theorem and also its Corollary.\\

\vskip 1pc \nind
{\bf Remark 3.7.} Now we are left with the case $\lkd(X)=-\infty$,
${\lkd}(B) = - \infty$ and
$\pi_1(X)$ is finite. This case appears to be much harder.\\

For example, let $X$ be the affine surface in $\C^3$ defined by
$\{XY-Z^2=1\}$. It is well-known that $X\cong\BP^1\times\BP^1-$ Diagonal. Clearly, $X$ has an
$\A^1$-fibration over $\BP^1$ which is Zariski-locally trivial. $X$ is
simply-connected and
$\lkd(X)=-\infty$. The fibration $X\rar\A^1$ given by the function $x$ is an $\A^1$
fibration
with
all reduced fibers. It is not clear even in this case if $X$
has a finite self-morphism of degree $>1$.\\

\newpage
\nind

\noindent
R. V. Gurjar, School of Mathematics, Tata Institute of Fundamental Research.
Homi-Bhabha road, Mumbai 400005, India.\\
e-mail: gurjar@math.tifr.res.in\\

\noindent
D. -Q. Zhang, Department of Mathematics, National University of Singapore, 2 Science Drive 2, Singapore
117543, Singapore.\\
e-mail: matzdq@nus.edu.sg\\

\end{document}